\documentclass[12pt]{amsart}
\usepackage{amsfonts}
\usepackage{amsmath}
\usepackage{amssymb}
\usepackage[margin=1in]{geometry}
\usepackage{tikz}
\usepackage{amsthm}   							
\usepackage{mathtools}							
\usepackage{commath}							
\usepackage{enumitem}  							

\usepackage[english]{babel}						

\makeatletter
\@namedef{subjclassname@2020}{\textup{2020} Mathematics Subject Classification}
\makeatother

\newtheorem{theorem}{Theorem}[section]

\newtheorem{lemma}[theorem]{Lemma}

\theoremstyle{definition}

\numberwithin{equation}{section}

\newcommand{\expm}{\exp^{\ast}}
\newcommand{\nconv}{^{\ast n}}

\newcommand{\N}{\mathbb{N}}

\renewcommand{\Re}{\operatorname{Re}}

\newcommand{\I}{\mathrm{i}}
\newcommand{\e}{\mathrm{e}}

\newcommand{\MP}{\mathcal{P}}
\newcommand{\MN}{\mathcal{N}}

\DeclareMathOperator{\Li}{Li}
\DeclareMathOperator*{\res}{Res}

\DeclareMathOperator{\supp}{supp}

\begin{document}

\title[Note on a conjecture of Bateman and Diamond concerning the abstract PNT]{Note on a conjecture of Bateman and Diamond concerning the abstract PNT with Malliavin-type remainder}

\author[F. Broucke]{Frederik Broucke}
\thanks{The author was supported by the Ghent University BOF-grant 01J04017}

\address{Department of Mathematics: Analysis, Logic and Discrete Mathematics\\ Ghent University\\ Krijgslaan 281\\ 9000 Gent\\ Belgium}
\email{fabrouck.broucke@ugent.be}


\begin{abstract}
Given $\beta\in(0,1)$, we show the existence of a Beurling generalized number system whose integer counting satisfies $N(x) = ax + O\bigl(x\exp(-c\log^{\beta} x)\bigr)$ for some $a>0$ and $c>0$, and whose prime counting function satisfies $\pi(x) = \Li(x) + \Omega\bigl(x\exp(-c'(\log x)^{\frac{\beta}{\beta+1}})\bigr)$ for some $c'>0$. This is done by generalizing a construction of Diamond, Montgomery, and Vorhauer. This Beurling system serves as additional motivation for a conjecture of Bateman and Diamond from 1969, concerning the PNT with Malliavin-type remainder.
\end{abstract}

\maketitle
\section{Introduction}
\label{sec: intro}
A.\ Beurling \cite{Beurling1937} defined a system of generalized primes $\MP$ as an unbounded sequence $(p_{j})_{j\ge1}$ satisfying $1< p_{1} \le p_{2} \le \dotso$. The corresponding system of generalized integers $\MN$ is the multiplicative semigroup generated by $1$ and $\MP$. With these systems one associates functions $\pi(x)$ and $N(x)$ counting the number of primes, respectively integers, not exceeding $x$, taking multiplicities into account for those integers having multiple representations as a product of generalized primes.

A substantial part of the theory of Beurling generalized primes focuses on the asymptotic behavior of the counting functions $\pi$ and $N$ and how one is related to the other. This was in fact Beurling's main motivation; in his seminal paper \cite{Beurling1937}, he showed the following abstract form of the prime number theorem (PNT). Suppose that for some $a>0$ and $\gamma>3/2$, $N(x) = ax + O\bigl(x\log^{-\gamma}x\bigr)$. Then $\pi(x)\sim x/\log x$ holds. In the same article, it was also shown that the condition $\gamma > 3/2$ is sharp (see also \cite{Diamond1970b}). The first abstract PNT with remainder actually predates Beurling's paper in some sense, since the methods of Landau's work on the prime ideal theorem \cite{Landau1903}, cast into the language of Beurling generalized primes, yield the following theorem (PNT with de la Vall\'ee Poussin remainder). Suppose that $N(x) = ax + O(x^{\theta})$ holds for some $a>0$ and $\theta<1$. Then there is a constant $c>0$ such that $\pi(x) = \Li(x) + O\bigl(x\exp(-c\sqrt{\log x})\bigr)$. Here $\Li$ represents the logarithmic integral, which we choose to define in this article as
\[
	\Li(x) \coloneqq \int_{1}^{x}\frac{1-u^{-1}}{\log u}\dif u.
\]
Landau's PNT was shown to be optimal by Diamond, Montgomery, and Vorhauer \cite{DiamondMontgomeryVorhauer} in the case $\theta>1/2$. Namely, they showed the existence of a system of Beurling generalized primes satisfying $N(x) = ax + O(x^{\theta})$ and $\pi(x) = \Li(x) + \Omega\bigl(x\exp(-c\sqrt{\log x})\bigr)$. (The relation $f = g + \Omega(h)$ means that $(f-g)/h \not\to 0$.) 

More general remainders (which we shall refer to as Malliavin remainders) were considered by Malliavin, namely
\begin{equation}
\label{eq: Malliavin PNT}
\tag{P$_{\alpha}$}
\pi(x)= \Li(x)+ O\bigl(x\exp (-c \log^{\alpha} x )\bigr) \qquad \mbox{for some } c > 0,
\end{equation}
and
\begin{equation}
\label{eq: Malliavin density}
\tag{N$_{\beta}$}
N(x)=  ax+ O\bigl(x\exp (-c' \log^{\beta} x )\bigr) \qquad \mbox{for some $a$, $c'>0$,} 
\end{equation}
with $0 < \alpha$, $\beta \le 1$.
Malliavin showed \cite{Malliavin} that \eqref{eq: Malliavin PNT} implies \eqref{eq: Malliavin density} with $\beta=\alpha/(\alpha+2)$, and that \eqref{eq: Malliavin density} implies \eqref{eq: Malliavin PNT} with $\alpha=\beta/10$. 
The first result was improved by Diamond \cite{Diamond1970a}, who showed that \eqref{eq: Malliavin PNT} implies \eqref{eq: Malliavin density} with $\beta=\alpha/(\alpha+1)$ and furthermore with $(\log x \log\log x)^{\beta}$ instead of $(\log x)^{\beta}$ in the exponential, see also \cite{HilberdinkLapidus2006} for the case $\alpha=1$. 
Diamond's improved result has recently been shown to be optimal for $\alpha=1$ by the author together with G.\ Debruyne and J.\ Vindas in \cite{B-D-V2020}. In an upcoming article, the same authors will show that the theorem of Diamond is optimal for every $0<\alpha\le 1$. 
In the converse direction, the value for $\alpha$ was improved to $\alpha = \beta/7.91$ by Hall in \cite{Hall1972}. 
(A slight refinement of his argument actually yields $\alpha=\beta/(\beta+6.91)$, se e.g.\ \cite[Section 16.4]{DiamondZhangbook}.) 
Hall's proof consists of a Tauberian argument combined with bounds on the zeta function which are obtained via a generalization of the familiar ``3-4-1-inequality''. The value $6.91$ arises from a specific choice of a positive trigonometric polynomial\footnote{The corresponding optimization problem for positive trigonometric polynomials dates back to Landau and is well-studied, see e.g.\ \cite{Revesz}. In particular it is known that the smallest value which can be obtained via this method is strictly above $6.87$.}.

Let us denote by $\beta^{\ast}(\alpha)$ the supremum of all $\beta$ for which \eqref{eq: Malliavin density} follows from \eqref{eq: Malliavin PNT}, and similarly $\alpha^{\ast}(\beta)$ as the supremum of all $\alpha$ for which \eqref{eq: Malliavin PNT} follows from \eqref{eq: Malliavin density}. The aforementioned results imply that $\beta^{\ast}(1) = 1/2$, $\beta^{\ast}(\alpha) \ge \alpha/(\alpha+1)$ for $0<\alpha<1$, $\alpha^{\ast}(1) = 1/2$ and $\alpha^{\ast}(\beta) \ge \beta/(\beta+6.91)$ for $0<\beta<1$. Based on Landau's PNT and the exponent $\beta=\alpha/(\alpha+1)$ in the converse problem, it was conjectured\footnote{They expressed this with the careful wording: ``There is quite likely room for improvement here (from $\alpha=\beta/7.91$), possibly to the value $\alpha=\beta/(\beta+1)$''.} by Bateman and Diamond \cite{BD69} that $\alpha^{\ast}(\beta) \ge \beta/(\beta+1)$.

In this note, we will show that $\alpha^{\ast}(\beta) \le \beta/(\beta+1)$. This will be done by showing the existence of a Beurling prime system satisfying
\[
	N(x) = ax + O\bigl(x\exp(-c'\log^{\beta} x)\bigr) \quad \mbox{and} \quad \pi(x) = \Li(x) + O\bigl(x\exp(-c(\log x)^{\frac{\beta}{\beta+1}})\bigr),
\]
for some $a>0$ and some $c$, $c'>0$, but for which also
\[
	\pi(x) = \Li(x) + \Omega\bigl(x\exp(-c''(\log x)^{\frac{\beta}{\beta+1}})\bigr), 
\]
for some $c''>c$, showing that the exponent $\beta/(\beta+1)$ cannot be improved. The example is found by generalizing the construction by Diamond, Montgomery, and Vorhauer, which corresponds to the case $\beta=1$. Since their example yields the optimal exponent $\alpha^{\ast}(1)=1/2$ in the case $\beta=1$, it is not unreasonable to imagine that the exponent $\beta/(\beta+1)$ occurring in a natural generalization of their example would also be optimal.

\section{Preliminaries}
\label{sec: preliminaries}
The construction of the example consists of two parts. First a \emph{continuous} generalized number system is provided by specifying its zeta function. This continuous system will satisfy the desired asymptotic relations. Then, using a probabilistic method due to Diamond, Montgomery, and Vorhauer, it is shown that there exists a \emph{discrete} system which approximates the continuous one in a suitable sense.

For a Beurling prime system $\MP$, one defines its zeta function $\zeta_{\MP}$ as $\zeta_{\MP}(s) = \sum_{n\in\MN}n^{-s}$, wherever this series converges. Here $s=\sigma + \I t$ is a complex variable. The function $\zeta_{\MP}$ satisfies the identity
\[
	\zeta_{\MP}(s) = \exp\biggl(\int_{1}^{\infty}x^{-s}\dif\Pi(x)\biggr),
\]
where $\Pi(x) \coloneqq \sum_{\nu\ge1}\pi(x^{1/\nu})/\nu$ (which is sometimes referred to as Riemann's prime counting function). The above identity is an immediate consequence of the defining property of $\MN$.

It is useful to extend the notion of generalized number system to include not necessarily discrete systems. In a broader sense \cite{Beurling1937}, \cite{DiamondZhangbook}, a Beurling generalized number system is defined as a pair $(\Pi, N)$ of non-decreasing, right-continuous, unbounded functions supported on $[1, \infty)$ and satisfying $\Pi(1)=0$, $N(1)=1$, and 
\[
	\zeta(s) \coloneqq \int_{1^{-}}^{\infty}x^{-s}\dif N(x) = \exp\biggl(\int_{1}^{\infty}x^{-s}\dif\Pi(x)\biggr).
\]
A system arising from a sequence $(p_{j})_{j\ge1}$, as defined in the introduction, will be called a \emph{discrete} generalized number system. When constructing examples of discrete systems, it is often easier to first construct a system (in the extended sense) with the desired properties, and to then ``discretize'' this system, rather than to come up with a discrete system right away. One possible way of discretizing systems is by a procedure used in \cite{DiamondMontgomeryVorhauer}, and later refined in \cite{Zhang2007}. This procedure yields the following
\begin{lemma}[Diamond, Montgomery, Vorhauer \cite{DiamondMontgomeryVorhauer}, Zhang \cite{Zhang2007}]
\label{lem: discretization}
	Let $f$ be a non-negative function supported on $[1,\infty)$ and satisfying
	\[
		f(u) \ll \frac{1}{\log u} \quad \mbox{and} \quad \int_{1}^{\infty}f(u)\dif u = \infty.
	\]
	Then there exists an increasing sequence of numbers $(p_{j})_{j\ge1}$, $p_{1}>1$ and $p_{j}\to\infty$ such that for any $t$ and any $x\ge1$
	\begin{equation}
	\label{eq: bound discretization}
		\abs{\sum_{p_{j}\le x}p_{j}^{-\I t} - \int_{1}^{x}u^{-\I t}f(u)\dif u } \ll \sqrt{x} + \sqrt{\frac{x\log(\,\abs{t}+1)}{\log(x+1)}}.
	\end{equation}
\end{lemma}
\noindent Note in particular that setting $t=0$ yields $\pi(x) = \sum_{p_{j}\le x}1 = \int_{1}^{x}f(u)\dif u + O(\sqrt{x})$.

\mbox{}

As said before, our construction is a natural generalization of the construction in \cite{DiamondMontgomeryVorhauer}. The results in \cite{DiamondMontgomeryVorhauer} were later sharpened by Zhang in \cite{Zhang2007}, and are also provided in the monograph \cite{DiamondZhangbook} of Diamond and Zhang. In this article, we shall roughly follow the structure of \cite[Sections 17.4-17.9]{DiamondZhangbook}. All necessary definitions and lemmas will be given; however if the proof of a statement is identical or very similar to a proof in \cite{DiamondZhangbook}, we will omit it and refer to \cite{DiamondZhangbook} instead. In fact, most of the arguments employed there can be carried over to the case $0<\beta<1$. However, a new difficulty arises as the considered zeta function does not seem to have meromorphic continuation beyond $\Re s=1$. We overcome this difficulty by considering natural approximations of it by more regular zeta functions.

\mbox{}

The main idea in \cite{DiamondMontgomeryVorhauer} is to construct a zeta function which has infinitely many zeros on the curve $\sigma=1-1/\log \abs{t}$, and none to the right of it, and which is of not too large growth. The zeros are ``responsible'' for the de la Vall\'ee Poussin remainder in the PNT, while the moderate growth of the zeta function allows one to deduce the desired asymptotics of $N$ via Perron inversion. We will modify the construction so that the zeros will lie on the curve\footnote{This curve is suggested by a result of Ingham, which shows what error term in the classical PNT would follow from a general zero-free region $\sigma>1-\eta(t)$ for the Riemann zeta function, see \cite[pages 60--65]{Inghambook}.} $\sigma=1-1/(\log\abs{t})^{1/\beta}$. The zeros are obtained by taking products of rescaled and translated versions of the function $G$ defined as
\[
	G(z) \coloneqq 1 - \frac{\e^{-z}-\e^{-2z}}{z}, \quad G(0)\coloneqq 0.
\]
The function $G$ has zeros $z_{0} = 0$ and $z_{\pm n} = x_{n} \pm \I y_{n}$, $n\in\N_{0}$, with 
\[
	-b\log \frac{n\pi}{2} \le x_{n} < -\frac{1}{2}\log\frac{n\pi}{2}, \quad \mbox{and} \quad n\pi < y_{n} < (n+1)\pi
\]
for some constant $b>1/2$. It has no other zeros. For $z\neq0$ we have the simple approximation
\begin{equation}
\label{eq: approx G}
	G(x+\I y) = 1+\theta\,\frac{\e^{-x}+\e^{-2x}}{\abs{x+\I y}}, \quad \abs{\theta} \le 1;
\end{equation}
while for $x \ge -1$, 
\begin{equation}
\label{eq: bound G}
\abs{G(x+\I y)} \le 1+\e^{2}-\e,
\end{equation}
which follows from the identity $G(z) = 1 - \int_{1}^{2}\e^{-zu}\dif u$. The logarithm of $G$ can also be expressed as a Mellin transform:
\begin{lemma}
\label{lem: log G}
The function $\log G(z)$ is well-defined for $x=\Re z>0$ and has the representation
\[
	\log G(z) = - \int_{1}^{\infty}g(u)u^{-z-1}\dif u,
\]
where
\[
	g(u) \coloneqq \sum_{n=1}^{\infty}\frac{1}{n}\chi\nconv(u).
\]
Here, $\chi$ is the indicator function of $[\e, \e^{2}]$ and $\ast$ denotes the multiplicative convolution of functions supported on $[1,\infty)$: $(f\ast h)(x) \coloneqq \int_{1}^{x}f(x/u)h(u)\dif u/u$. The function $g$ is non-negative, supported on $[\e, \infty)$, and on intervals $(\e^{m},\e^{m+1})$ it equals a polynomial in $\log u$ of degree at most $m-1$.
\end{lemma}
The function $g(u)$ gets close to $1/\log u$ for large $u$. We have the following estimates.

\begin{lemma}
\label{lem: estimate g}
	For $u>\e^{2}$, 
	\[	g(u)\log u = 1 + O\bigl(u^{-(1/2)\log(\pi/2)}\bigr), \]
	and for $u\ge\e^{5}$, $g$ is differentiable and satisfies
	\[	(g(u)\log u)' = O\bigl(u^{-1-(1/2)\log(\pi/2)}\bigr). \]
\end{lemma}
\noindent For proofs of the above statements and lemmas, we refer to \cite[Sections 17.5, 17.6]{DiamondZhangbook}.

\section{The example}
\label{sec: the example}
We will now define the continuous example by specifying its zeta function. Let $\beta\in(0,1)$ and set
\[
	l_{k}=4^{k}, \quad \gamma_{k}=\exp\bigl((l_{k})^{\beta}\bigr) = \e^{4^{\beta k}}, \quad \rho_{k}=1-\frac{1}{l_{k}}+\I\gamma_{k}.
\]
These are the same parameters as in \cite{DiamondZhangbook}, except for $\gamma_{k}$, which we have set to be $\exp\bigl((l_{k})^{\beta}\bigr)$ instead of $\exp(l_{k})$. The points $\rho_{k}$ now lie on the curve $\sigma=1-1/(\log t)^{1/\beta}$ instead of $\sigma=1-1/\log t$. Next we set
\[
	\zeta_{C}(s) \coloneqq \frac{s}{s-1}\prod_{k=1}^{\infty}G(l_{k}(s-\rho_{k}))G(l_{k}(s-\overline{\rho}_{k})).
\]
Using \eqref{eq: approx G}, we see that the product converges uniformly in the half plane $\sigma\ge1$, so this zeta function is holomorphic in the open half plane $\sigma>1$. For $\beta<1$, this zeta function does not seem to have analytic continuation to a larger half plane, unlike the case $\beta=1$. The factor $s/(s-1)$ corresponds to the main term $\Li(x)$ in the PNT, while the factors of the infinite product will produce the desired oscillation. That $\zeta_{C}$ is indeed the zeta function of a Beurling system follows from the following lemma.
\begin{lemma}
\label{lem: def f_{C}}
	For $\sigma>1$, 
	\[	\zeta_{C}(s) = \exp\biggl(\int_{1}^{\infty}v^{-s}f_{C}(v)\dif v\biggr), \]
	with
	\begin{equation}
	\label{eq: def f_{C}}
	f_{C}(v) \coloneqq \frac{1-v^{-1}}{\log v} - 2\sum_{k\ge1}\frac{g(v^{1/l_{k}})}{l_{k}}v^{-1/l_{k}}\cos(\gamma_{k}\log v), \quad v\ge1.
	\end{equation}
	We have $f_{C}(v)>0$ for $v>1$ and $f_{C}$ satisfies the Chebyshev estimates for some $\delta\in(0,1)$
	\[	(1-\delta)\frac{1-v^{-1}}{\log v} \le f_{C}(v) \le (1+\delta)\frac{1-v^{-1}}{\log v}, \quad v\ge\e^{4}. \]	
\end{lemma}
\noindent The proof is identical to the one in \cite{DiamondZhangbook}, since the cosine is bounded trivially by $1$, and this is the only place where the altered parameter $\gamma_{k}$ occurs.

The lemma implies that $\zeta_{C}$ is the zeta function of the Beurling number system with Riemann prime counting function $\Pi_{C}$ given by $\Pi_{C}(x)=\int_{1}^{x}f_{C}(v)\dif v$. The integer counting function $N_{C}$ is uniquely determined by $\Pi_{C}$ (explicitly $\dif N_{C} = \expm(\dif\Pi_{C})$, where the exponential is the exponential with respect to multiplicative convolution of measures on $[1, \infty)$ (see e.g. \cite[Chapter 3]{DiamondZhangbook})).

Next, we use Lemma \ref{lem: discretization} with $f=f_{C}$ to obtain a sequence $\MP = (p_{j})_{j\ge1}$ of Beurling primes which satisfies \eqref{eq: bound discretization} with $f=f_{C}$. Denote the prime and integer counting function of $\MP$ by $\pi$ and $N$ respectively. We also consider its Chebyshev prime counting function $\psi$, defined as the summatory function of $\Lambda$, with $\Lambda(n)=\log p_{j}$ if $n=p_{j}^{\nu}$ for some $p_{j}\in\MP$, $\nu\ge1$, and zero otherwise. In the next two sections, we will show the following relations:
\begin{equation}
	N(x) = ax + O\bigl(x\exp(-c\log^{\beta} x)\bigr), 	\quad \mbox{for some $a>0$ and $c>0$}; \label{eq: asymptotics N} 
\end{equation}
and
\begin{equation}
	\begin{split}
		&\limsup_{x\to\infty}\frac{\psi(x)-x}{x\exp\bigl(-\beta^{-\frac{\beta}{\beta+1}}(\beta+1)(\log x)^{\frac{\beta}{\beta+1}}\bigr)} = 2,\\
		&\liminf_{x\to\infty}\frac{\psi(x)-x}{x\exp\bigl(-\beta^{-\frac{\beta}{\beta+1}}(\beta+1)(\log x)^{\frac{\beta}{\beta+1}}\bigr)} = -2.
	\end{split}  \label{eq: asymptotics psi}
\end{equation}
From these two relations it then follows that $\alpha^{\ast}(\beta) \le \beta/(\beta+1)$, since $\pi(x) = \int_{1}^{x}(1/\log u) \dif\psi(u) + O(\sqrt{x})$.

\section{Asymptotics of $N$}
\label{sec: asymptotics N}
In order to deduce asymptotic information of $N$, we will use a Perron inversion formula. We will bypass the problem of the apparent absence of analytic continuation of $\zeta_{C}$ beyond $\sigma=1$ by considering for $K\ge1$
\[
	\zeta_{C,K}(s) \coloneqq \frac{s}{s-1}\prod_{k=1}^{K}G(l_{k}(s-\rho_{k}))G(l_{k}(s-\overline{\rho}_{k})) = \exp\biggl(\int_{1}^{\infty}v^{-s}f_{C,K}(v)\dif v\biggr),
\]
where $f_{C,K}$ is defined as in \eqref{eq: def f_{C}}, but with the summation ranging only up to $K$. Then $\zeta_{C,K}$ has meromorphic continuation to the whole complex plane. Note that also $f_{C,K}>0$, and that $f_{C}(v) = f_{C,K}(v)$ whenever $v < \e^{4^{K+1}}$, since $\supp g(v^{1/l_{k}}) \subseteq [\e^{l_{k}}, \infty)$. Next we set
\begin{equation}
\label{eq: def Pi_{K}}	
	\dif \Pi_{K}(v) = \chi_{[1,\e^{4^{K+1}})}(v)\dif\Pi(v) + \chi_{[\e^{4^{K+1}}, \infty)}(v)f_{C,K}(v)\dif v. 
\end{equation}
Here, $\Pi$ is the Riemann prime counting function of the discrete system $\MP$, and $\chi_{I}$ denotes the indicator function of a set $I$. With $\Pi_{K}$ we associate the integer counting function $N_{K}$ (i.e. $\dif N_{K} = \expm(\dif\Pi_{K})$) and the zeta function $\zeta_{K}$. One might view the Beurling system $(\Pi_{K},N_{K})$ as an intermediate system between the discrete system $\MP$ and the continuous one given by $f_{C,K}$. Since $\Pi_{K} = \Pi$ on $[1,\e^{4^{K+1}})$, also $N_{K}=N$ on $[1,\e^{4^{K+1}})$. 

We will apply the following Perron formula for $\int N_{K}$ (to guarantee absolute convergence):
\begin{equation}
\label{eq: Perron formula}
	\int_{1}^{x}N_{K}(u)\dif u = \frac{1}{2\pi\I}\int_{\kappa-\I\infty}^{\kappa+\I\infty}\zeta_{K}(s)x^{s+1}\frac{\dif s}{s(s+1)}.
\end{equation}
Here $\kappa$ is a number larger than $1$. We will shift the contour to a contour to the left of $\sigma=1$, picking up a residue at $s=1$ which will provide the main term in \eqref{eq: asymptotics N}. The integral over the shifted contour will be estimated by comparing $\zeta_{K}$ with $\zeta_{C,K}$, and by applying some bounds for $\zeta_{C,K}$ which we will derive shortly. We have
\begin{align*}
	\log \zeta_{K}(s) - \log\zeta_{C,K}(s) 		&= \int_{1}^{\e^{4^{K+1}}}u^{-s}(\dif\Pi_{K}(u) - f_{C,K}(u)\dif u)  \\
									&= \int_{1}^{\e^{4^{K+1}}}u^{-s}(\dif\Pi(u)-\dif\pi(u)) + \int_{1}^{\e^{4^{K+1}}}u^{-\sigma}u^{-\I t}(\dif\pi(u) - f_{C}(u)\dif u).
\end{align*}
Here we used that $f_{C,K} = f_{C}$ on $[1,\e^{4^{K+1}})$. Note that both of these integrals are entire functions of $s$. Since $\dif\Pi-\dif\pi$ is a positive measure, and since $\Pi(x)-\pi(x) = O(\sqrt{x})$ (this follows immediately from the definition of $\Pi$ and the bound $\pi(x) \ll x$), the first integral is uniformly bounded (independent of $K$) in the half plane $\sigma\ge3/4$ say. For the second integral, we integrate by parts and use the bound \eqref{eq: bound discretization} to see that it is uniformly bounded by a constant times $\sqrt{\log(2+\abs{t})}$ in the half plane $\sigma\ge3/4$. For the remainder of this section we fix positive constants $A$ and $B$, independent of $K$, so that
\begin{equation}
\label{eq: bound log zetas}
	\abs{\log \zeta_{K}(s) - \log\zeta_{C,K}(s)} \le \begin{cases}
						A 					&\mbox{ if $\sigma\ge3/4$, $\abs{t}\le 2$} \\
						A + B\sqrt{\log\abs{t}} 	&\mbox{ if $\sigma\ge3/4$, $\abs{t}\ge2$.}  
						\end{cases}
\end{equation}

Let now $x\ge\e^{4}$ be fixed, and let $K$ be such that $\e^{4^{K}} \le x < \e^{4^{K+1}}$. Then $N(x)=N_{K}(x)$. Set $\sigma_{1}=1-(1/2)(\log x)^{\beta-1}$, $\sigma(t) = 1-(1/4)\log \abs{t}/\log x$, and let $k(\beta)$ be such that $(3/2)\gamma_{k} < (1/2)\gamma_{k+1}$ for $k\ge k(\beta)$.
\begin{lemma}
\label{lem: bound zeta_{K}}
The following bounds hold uniformly (with implicit constants independent of $K$):

\noindent for $\sigma_{1}\le\sigma\le2$:
\begin{enumerate}
	\item if $0\le t \le 2$, then $\zeta_{K}(\sigma+\I t) \ll 1/\abs{\sigma-1}$;
	\item if $t\ge 2$ and $\abs{t-\gamma_{k}} \ge \gamma_{k}/2$, for every $k\in \{k(\beta), k(\beta)+1, \dotso, K\}$, then 
	\[	\zeta_{K}(\sigma+\I t) \ll \exp\bigl(B\sqrt{\log\abs{t}}\bigr); \]
	\item if $t\ge2$ and $\abs{t-\gamma_{k_{0}}} < \gamma_{k_{0}}/2$ for some $k_{0}\in\{ k(\beta), k(\beta)+1,\dotso,K\}$, then 
	\[
		\zeta_{K}(\sigma+\I t) \ll \begin{dcases}
						\exp\bigl(B\sqrt{\log\abs{t}}\bigr)\biggl(1+\frac{\abs{t}}{4^{k_{0}}\abs{\sigma+\I t-\rho_{k_{0}}}}\biggr)	&\mbox{if } 4^{k_{0}}\abs{\sigma+\I t-\rho_{k_{0}}} \ge 1, \\
						\exp\bigl(B\sqrt{\log\abs{t}}\bigr)													&\mbox{if } 4^{k_{0}}\abs{\sigma+\I t-\rho_{k_{0}}} \le 1;
						\end{dcases}
	\]
\end{enumerate}
for $\abs{t}\ge 2\gamma_{K}$ and $\sigma(t) \le \sigma \le 2$:
\[	\zeta_{K}(\sigma+\I t) \ll \exp\bigl(B\sqrt{\log\abs{t}}\bigr). \]
\end{lemma}
The proof is essentially the same as that of \cite[Lemma 17.22]{DiamondZhangbook}. First we use \eqref{eq: bound log zetas} to compare with $\zeta_{C,K}$. Then we use \eqref{eq: approx G} to approximate the factors in the product. The main point is that $\exp(2l_{k}(1-\sigma)) \ll \gamma_{k}$ for $k\le K$ for $\sigma\ge\sigma_{1}$, while for $\sigma\ge\sigma(t)$, we have $\exp(2l_{k}(1-\sigma)) \le \sqrt{\,\abs{t}}$. By definition of $k(\beta)$, for each fixed $t$ there is at most one $k_{0}\in\{k(\beta), \dotso,K\}$ with $\abs{t-\gamma_{k_{0}}}<\gamma_{k_{0}}/2$. For the terms with $k<k(\beta)$, we just employ some uniform bound in the half plane $\sigma\ge3/4$ say. We omit the details.

\mbox{}

Let us now focus our attention on the Perron integral. Since $N_{K}$ is non-decreasing, $\int_{x-1}^{x}N_{K}(u)\dif u \le N_{K}(x) \le \int_{x}^{x+1}N_{K}(u)\dif u$. Combining this with the Perron formula \eqref{eq: Perron formula}, we have
\[
	N_{K}(x) \le \frac{1}{2\pi\I}\int_{\kappa-\I\infty}^{\kappa+\I\infty}\zeta_{K}(s)\frac{(x+1)^{s+1}-x^{s+1}}{s(s+1)}\dif s.
\]
We shift the contour to a contour $\Gamma$ to the left of $\sigma =1$. Set $J=\lfloor (\log x)^{\min(3\beta/2, 1)} \rfloor$ and 
\begin{align*}
	\Gamma_{1} 	&\coloneqq [\sigma_{1}, \sigma_{1}+\I2^{J}]; \\ 
	\Gamma_{2} 	&\coloneqq [\sigma_{1}+\I2^{J}, \sigma(2^{J})+\I2^{J}] \cup \{\sigma(t)+\I t: 2^{J}\le t\le x\}; \\	
	\Gamma_{3} 	&\coloneqq [3/4 + \I x, 3/4 + \I\infty).
\end{align*}
Note that $\sigma_{1}>\sigma(2^{J})$. We let $\Gamma$ be the union of the $\Gamma_{i}$ and their complex conjugates. Moving from $\sigma=\kappa$ to $\Gamma$ is allowed, since the contribution of the connecting piece $[3/4+\I T, \kappa+\I T]$ tends to 0 as $T\to\infty$, by the last bound of Lemma \ref{lem: bound zeta_{K}}. By the residue theorem we have 
\[
	N_{K}(x) \le a_{K}(x+1/2) + O(I_{1}+I_{2}+I_{3}),
\]
where $a_{K} = \res_{s=1}\zeta_{K}(s)$ and $I_{i}$ the integral over $\Gamma_{i}$. For $I_{1}$ we perform a dyadic splitting, and write
\[
	I_{1} = I_{1,0} + \sum_{j=1}^{J-1}I_{1,j}, 
\]
where $I_{1,0}$ is the part of the integral where $0\le t\le 2$, and $I_{1,j}$ the part of the integral where $2^{j}\le t\le 2^{j+1}$. By bounding $(x+1)^{s+1}-x^{s+1}$ by $\abs{s+1}x^{\sigma}$, and using the bounds from Lemma \ref{lem: bound zeta_{K}}, we see that
\[
	I_{1,0} \ll x\exp\bigl(-(1/2)\log^{\beta} x\bigr)(\log x)^{1-\beta}, \quad I_{1,j} \ll x\exp\bigl(-(1/2)\log^{\beta} x\bigr)\exp\bigl(B\sqrt{\log 2^{j+1}}\bigr).
\] 
Indeed, if some $t$ satisfies $\abs{t-\gamma_{k_{0}}} < \gamma_{k_{0}}/2$ for some $k_{0}\ge k(\beta)$, we use the third bound of Lemma \ref{lem: bound zeta_{K}} and the estimate
 \[
 	\int\limits_{\substack{\abs[0]{t-\gamma_{k_{0}}} < \gamma_{k_{0}}/2, \\ \abs[0]{\sigma_{1}+\I t -\rho_{k_{0}}} \ge 4^{-k_{0}}}} \frac{\dif t}{4^{k_{0}}\abs{\sigma_{1}+\I t -\rho_{k_{0}}}} 
	\ll \frac{1}{4^{k_{0}}}\bigl(\log \gamma_{k_{0}} + \log 4^{k_{0}}\bigr) \ll 1.
\]
Also, by definition of $k(\beta)$, there are at most two values of $k\ge k(\beta)$ such that $\abs{t-\gamma_{k}}<\gamma_{k}/2$ for some $t$ in a dyadic interval $[2^{j},2^{j+1}]$.
Using the estimate $\sum_{j\le J}\e^{D\sqrt{j}} \ll \sqrt{J}\e^{D\sqrt{J}}$ we get
\[
	I_{1} \ll x\exp\bigl(-(1/2)\log^{\beta} x\bigr)(\log x)\exp\bigl(O\bigl((\log x)^{3\beta/4}\bigr)\bigr) \ll x\exp\bigl(-c\log^{\beta} x\bigr),
\]
for any $c<1/2$.

For $I_{2}$, we again bound $(x+1)^{s+1}-x^{s+1}$ by $\abs{s+1}x^{\sigma}$ and use the last bound of Lemma \ref{lem: bound zeta_{K}} (note that $2^{J} \ge 2\gamma_{K}$). We get
\begin{align*}
	I_{2} 	&\ll x\exp\bigl(-(1/2)\log^{\beta}x\bigr) + x\int_{2^{J}}^{x}\exp\bigl(B\sqrt{\log t} - (5/4)\log t\bigr)\dif t \\
		&\ll x\exp\bigl(-(1/2)\log^{\beta}x\bigr) + x\exp\bigl(-((\log2)/8)(\log x)^{\min(3\beta/2,1)}\bigr) \ll x\exp\bigl(-(1/2)\log^{\beta}x\bigr).
\end{align*}
Lastly, we bound $(x+1)^{s+1}-x^{s+1}$ by $x^{\sigma+1}$ and use the last bound from Lemma \ref{lem: bound zeta_{K}} to get
\[
	I_{3} \ll x^{7/4}\int_{x}^{\infty}\exp\bigl(B\sqrt{\log t}-2\log t\bigr)\dif t \ll x^{7/8}.
\]
Concluding the above calculations, we have that for any $c<1/2$
\[
	N_{K}(x) \le a_{K}x + O\bigl(x\exp(-c\log^{\beta}x)\bigr).
\]
The inequality $\ge$ can be shown in a completely analogous way. 

It remains to see that $a_{K}$ is close to $a$, the density of $N$. We have that
\[
	a_{K} = \exp\biggl(\int_{1}^{\infty}\frac{1}{u}\bigl(\dif\Pi_{K}(u)-\dif\Li(u)\bigr)\biggr), \quad a = \exp\biggl(\int_{1}^{\infty}\frac{1}{u}\big(\dif\Pi(u)-\dif\Li(u)\bigr)\biggr),
\]
where we used that $\exp\int_{1}^{\infty}u^{-s}\dif\Li(u) = s/(s-1)$ to compute the residues\footnote{If Beurling generalized integers have a density, it is equal to the right hand residue of its zeta function, see e.g.\ \cite[Proposition 5.1]{DiamondZhangbook}.}. The fact that both integrals converge follows from the fact that $\Pi_{K}(x) - \Li(x)$, $\Pi(x) - \Li(x) \ll x\exp\bigl(-c'\log^{\alpha} x\bigr)$ for some $\alpha>0$ and $c'>0$, cfr. Section \ref{sec: asymptotics of psi}.
By \eqref{eq: def Pi_{K}},
\begin{align*}
	a_{K} 	&= a\exp\biggl(\int_{\e^{4^{K+1}}}^{\infty}\frac{1}{u}\bigl(f_{C,K}(u)\dif u - \dif\Pi(u)\bigr)\biggr)\\
			&= a\exp\biggl(\int_{\e^{4^{K+1}}}^{\infty}\frac{1}{u}\bigl(f_{C,K}(u)-f_{C}(u)\bigr)\dif u + \int_{\e^{4^{K+1}}}^{\infty}\frac{1}{u}\bigl(f_{C}(u)\dif u - \dif\Pi(u)\bigr)\biggr).
\end{align*}
The first integral equals 
\[
	-\sum_{k=K+1}^{\infty} \log \bigl( G(l_{k}(1-\rho_{k}))G(l_{k}(1-\overline{\rho}_{k}))\bigr) \ll \sum_{k=K+1}^{\infty}\frac{1}{l_{k}\gamma_{k}} \ll \e^{-4^{\beta(K+1)}} \ll \exp\bigl(-\log^{\beta}x\bigr),
\]
where we used \eqref{eq: approx G}, and the second integral is bounded by $1/\sqrt{\e^{4^{K+1}}} \ll 1/\sqrt{x}$ by \eqref{eq: bound discretization}.
This gives that $a_{K} = a\bigl\{1+O\bigl(\exp(-\log^{\beta} x)\bigr)\bigr\}$. We conclude that for any $c<1/2$, 
\[
	N(x) = N_{K}(x) = a_{K}x + O\bigl(x\exp(-c\log^{\beta}x)\bigr) = ax + O\bigl(x\exp(-c\log^{\beta}x)\bigr),
\]
which shows \eqref{eq: asymptotics N}.

\section{Asymptotics of $\psi$}
\label{sec: asymptotics of psi}
The analysis of the prime counting function of the Diamond-Montgomory-Vorhauer example (corresponding to $\beta=1$) can be readily adapted to the case of general $\beta$; no new technical difficulties arise. We give a summary of the analysis, but refer to \cite[Section 17.9]{DiamondZhangbook} for the details.

Given a fixed $x$, let $K$ again be such that $\e^{4^{K}}\le x < \e^{4^{K+1}}$. We shall analyze the Chebyshev prime counting function $\psi(x) = \int_{1}^{x}\log u\dif \Pi(u)$. First note that $\psi(x) = \psi_{C}(x) + O(\sqrt{x}\log x)$, which follows from Lemma \ref{lem: discretization}, so it suffices to analyze $\psi_{C}$. Using the same notations as in \cite{DiamondZhangbook}, we have
\[
	\psi_{C}(x) = x -1 - \log x - 2F(x),
\]
with 
\[
	F(x) = \sum_{k=1}^{K}\int_{\e^{4^{k}}}^{x}(\log v)4^{-k}g(v^{4^{-k}})v^{-4^{-k}}\cos(\gamma_{k}\log v)\dif v \eqqcolon \sum_{k=1}^{K}I_{k}(x).
\]
Transforming the integrals with the substitution $u=v^{4^{-k}}$, splitting the integration range in $[\e, \e^{5}]$ and $[\e^{5}, x^{4^{-k}}]$, integrating by parts, and using the bounds from Lemma \ref{lem: estimate g}, one shows that 
\[
	I_{k}(x) = \frac{x^{1-4^{-k}}}{\gamma_{k}}\sin(\gamma_{k}\log x) + O\biggl\{ x^{5/16} + 
		\frac{1}{\gamma_{k}}\biggl(x^{1-4^{-k}(1+\frac{1}{2}\log(\pi/2))} + \frac{x^{1-4^{-k}}}{\gamma_{k}}\biggr)\biggr\}, \quad k\le K-2.
\]
To estimate the integrals $I_{K-1}$, $I_{K}$, we transform again to the variable $u=v^{4^{-k}}$ and split the integration range in intervals $[\e^{m}, \e^{m+1})$, $m < 16$. On $[\e^{m}, \e^{m+1})$, we write $g(u)$ as a polynomial in $\log u$ of degree at most $m-1$, and integrate by parts. This yields
\[
	I_{K-1}(x) \ll x\exp(-4^{-2\beta}\log^{\beta}x), \quad I_{K}(x) \ll x\exp(-4^{-\beta}\log^{\beta}x),
\]
which is ok, since $\beta>\beta/(\beta+1)$.
One then proceeds by showing that $F(x)$ is dominated by at most two terms $I_{k_{0}}(x)$ and $I_{k_{0}+1}(x)$, with $k_{0}$ close to $\frac{1}{\beta+1}(K+\log(1/\beta))$. Consider
\[
	\frac{x^{1-4^{-k}}}{\gamma_{k}} = x\exp\bigl(-4^{-k}\log x - 4^{\beta k}\bigr) = x \exp\biggl(-\frac{\log x}{\lambda} - \lambda^{\beta}\biggr),
\]
where we have written $\lambda=4^{k}$. The function $-\lambda^{-1}\log x - \lambda^{\beta}$ reaches it maximum at $\lambda_{\mathrm{max}}$,
\[
	\lambda_{\mathrm{max}} = \biggl(\frac{\log x}{\beta}\biggr)^{\frac{1}{\beta+1}}, \quad \mbox{and } 
	-\frac{\log x}{\lambda_{\mathrm{max}}} - (\lambda_{\mathrm{max}})^{\beta} = - \beta^{-\frac{\beta}{\beta+1}}(\beta+1)(\log x)^{\frac{\beta}{\beta+1}}.
\]
Note that $\lambda_{\mathrm{max}} < 4^{K-2}$ for $x$ sufficiently large. Now set $\mu = \log\lambda_{\mathrm{max}}/\log 4$, $k_{0} = \lfloor \mu \rfloor$, and write 
\[
	E(x) = x\exp\bigl(-\beta^{-\frac{\beta}{\beta+1}}(\beta+1)(\log x)^{\frac{\beta}{\beta+1}}\bigr).
\]
We have
\begin{align*}
	&\abs{I_{k_{0}}(x)} \le E(x)\{1+o(1)\},  	&&I_{k_{0}+1}(x) = o(E(x)) 				&& \mbox{if }\{\mu\} \in [0,1/3]; \\
	&I_{k_{0}}(x) = o(E(x)),				&&I_{k_{0}+1}(x) = o(E(x))					&& \mbox{if }\{\mu\} \in (1/3, 2/3); \\
	&I_{k_{0}}(x) = o(E(x)),  				&&\abs{I_{k_{0}+1}(x)} \le E(x)\{1+o(1)\}		&& \mbox{if }\{\mu\} \in [2/3, 1).		
\end{align*}
Also in every case, the terms $I_{k}(x)$, $k\neq k_{0}, k_{0}+1$ are $O\bigl(x\exp\bigl(-d(\log x)^{\frac{\beta}{\beta+1}}\bigr)\bigr)$ for some $d>\beta^{-\frac{\beta}{\beta+1}}(\beta+1)$, and there are $K-2 = O(\log\log x)$ such terms. Combining all these estimates shows that 
\[
	\limsup_{x\to \infty} \frac{\psi(x)-x}{E(x)} \le 2, \qquad \liminf_{x\to \infty} \frac{\psi(x)-x}{E(x)} \ge -2.
\]
In order to show equality and hence prove \eqref{eq: asymptotics psi}, one considers an increasing sequence of values for $x$, so that $(\log x/\beta)^{\frac{1}{\beta+1}}$ gets arbitrarily close to perfect fourth powers $4^{k_{0}}$, for some $k_{0} \le K-2$ ($k_{0}$ and $K$ of course depending on $x$), and where also $\sin(\gamma_{k_{0}}\log x)$ gets arbitrarily close to $-1$ (for the $\limsup$) or $1$ (for the $\liminf$).

\mbox{}

The author thanks dr.\ Gregory Debruyne and prof.\ Jasson Vindas for their kind advise and useful remarks.

\end{document}